# Accelerating Nuclear Configuration Interaction Calculations through a Preconditioned Block Iterative Eigensolver

Meiyue Shao[1], H. Metin Aktulga[2], Chao Yang[1], Esmond G. Ng[1], Pieter Maris[3], and James P. Vary[3]

[1]Computational Research Division, Lawrence Berkeley National Laboratory, Berkeley, CA 94720, United States
[2]Department of Computer Science and Engineering, Michigan State University, East Lansing, MI 48824, United States
[3]Department of Physics and Astronomy, Iowa State University, Ames, IA 50011, United States

August 25, 2017

## Abstract

We describe a number of recently developed techniques for improving the performance of large-scale nuclear configuration interaction calculations on high performance parallel computers. We show the benefit of using a preconditioned block iterative method to replace the Lanczos algorithm that has traditionally been used to perform this type of computation. The rapid convergence of the block iterative method is achieved by a proper choice of starting guesses of the eigenvectors and the construction of an effective preconditioner. These acceleration techniques take advantage of special structure of the nuclear configuration interaction problem which we discuss in detail. The use of a block method also allows us to improve the concurrency of the computation, and take advantage of the memory hierarchy of modern microprocessors to increase the arithmetic intensity of the computation relative to data movement. We also discuss implementation details that are critical to achieving high performance on massively parallel multi-core supercomputers, and demonstrate that the new block iterative solver is two to three times faster than the Lanczos based algorithm for problems of moderate sizes on a Cray XC30 system.

## 1 Introduction

The structure of an atomic nucleus with $A$ nucleons can be elucidated by solutions of the many-body Schrödinger equation

$$H\Psi(\mathbf{r}_1, \mathbf{r}_2, \ldots, \mathbf{r}_A) = \lambda \Psi(\mathbf{r}_1, \mathbf{r}_2, \ldots, \mathbf{r}_A), \qquad (1)$$

where $H$ is the nuclear Hamiltonian acting on the $A$-body wavefunction $\Psi(\mathbf{r}_1, \mathbf{r}_2, \ldots, \mathbf{r}_A)$ (with $\mathbf{r}_j \in \mathbb{R}^3$, $j = 1, 2, \ldots, A$), and $\lambda$ the energy. For stable nuclei, the low-lying spectrum is discrete. The solution associated with the lowest eigenvalue is the ground state. The nuclear Hamiltonian $H$ contains the kinetic energy operator $K$, and the potential term $V$ which describes the strong interactions between nucleons as well as the Coulomb interactions between protons.

Nuclear structure presents a unique set of challenges for achieving robust simulations within quantum many-body theory. First, the underlying strong interactions span scales from low energies (few MeV) to high energies (hundreds of MeV) as observed in the measured nucleon–nucleon scattering cross section. Second, many phenomena of interest, such as halo nuclei, collective rotation, fission and radioactive decay, involve the delicate interplay of binding contributions across those energy scales leading to these important low-energy and long-range phenomena. Third, in order to achieve accurate results from nuclear theory that are competitive with accurate experimental data, it is now recognized that we require three-nucleon (3N) interactions whose contributions represent one to two orders of magnitude increase in computational burden over calculations with two-nucleon (NN) interactions alone.

One approach to obtaining approximate solutions of (1) is to use the *No-Core Configuration Interaction* (NCCI) method [7] which projects the nuclear many-body Hamiltonian operator in (1) onto a finite-dimensional subspace spanned by *Slater determinants* of the form

$$\Phi_a(\mathbf{r}_1, \mathbf{r}_2, \ldots, \mathbf{r}_A) = \frac{1}{\sqrt{A!}} \det \begin{bmatrix} \phi_{a_1}(\mathbf{r}_1) & \phi_{a_2}(\mathbf{r}_1) & \ldots & \phi_{a_A}(\mathbf{r}_1) \\ \phi_{a_1}(\mathbf{r}_2) & \phi_{a_2}(\mathbf{r}_2) & \ldots & \phi_{a_A}(\mathbf{r}_2) \\ \vdots & \vdots & \ddots & \vdots \\ \phi_{a_1}(\mathbf{r}_A) & \phi_{a_2}(\mathbf{r}_A) & \ldots & \phi_{a_A}(\mathbf{r}_A) \end{bmatrix}, \qquad (2)$$

where $\phi_a$'s are orthonormal single-particle wavefunctions indexed by a generic label $a$. (Quantum mechanics dictates that wavefunctions of $A$ identical spin-$\frac{1}{2}$ particles are completely anti-symmetric, a condition satisfied by a Slater determinant. In practice, we often form our nuclear many-body wavefunctions as products of separate Slater determinants for $Z$ protons and $N$ neutrons, with $A = Z + N$.) [1]

The lowest eigenvalues and associated eigenvectors of the finite dimensional Hamiltonian $\hat{H}$, which is large and sparse, are then computed. The dimension of $\hat{H}$ depends on (1) the number of nucleons $A$; (2) the number of single-particle states; and (3) the (optional) many-body truncation. In the NCCI approach, typically one works in a basis of harmonic oscillator single-particle states where the number of single-particle states is implicitly determined by the many-body truncation $N_{\max}$, which imposes a limit on the sum of the single-particle energies (oscillator quanta) included in each Slater determinant of $A$ nucleons. In the limit of a complete (but infinite-dimensional) basis, this approach would give the exact bound state wave functions; in practice, increasingly accurate approximations to both the ground state and the narrow (low-lying) excited states of a given nucleus often require increasing $N_{\max}$ to the largest value where the eigenvalue problem is still solvable with available methods and resources. The dimension of $\hat{H}$ increases rapidly both with the number of nucleons, $A$, and with the truncation parameter, $N_{\max}$. This large sparse matrix eigenvalue problem has been demonstrated to be a computationally challenging problem on high-performance computers.

An efficient way to construct and diagonalize the projected Hamiltonian $\hat{H}$ has been implemented in the software package MFDn (Many Fermion Dynamics for nuclear structure) developed by Vary *et al.* [22, 23, 36]. MFDn uses the Lanczos algorithm [21, 28] to compute the desired eigenpairs. Over the last several years, we have developed a number of techniques to improve the computational efficiency of MFDn. These techniques include

1. an efficient scalable parallel scheme for constructing the matrix Hamiltonian [35];

2. efficient data distribution schemes that take into account the topology of the interconnect [3];

---

[1] Each Slater determinant corresponds to a quantum many-body configuration. The configurations interact with each other through the Hamiltonian (which contains both the kinetic and potential energy contributions), hence the name "Configuration Interaction". Since we include all nucleons of the nucleus in these configurations, we use the term "No-Core".

3. a technique to overlap communication with computation in a hybrid MPI/OpenMP programming model [4];

4. an efficient scheme to multiply the sparse matrix Hamiltonian with a number of vectors [2].

The last technique is extremely useful in block iterative algorithms for computing eigenvalues and eigenvectors of a sparse matrix. However, additional techniques are required in order to develop an efficient block eigensolver.

In this paper, we present implementation details of a preconditioned block iterative eigensolver to accelerate the NCCI computations in MFDn. Specifically, we adopt the Locally Optimal Block Preconditioned Gradient (LOBPCG) algorithm [17] as our eigensolver. The LOBPCG algorithm and its variants have been used to solve eigenvalue problems arising from a number of applications, including material science [6, 15, 20, 37] and machine learning [18, 19, 24]. The use of a block iterative method allows us to improve the memory access pattern of the computation and make use of approximations to several eigenvectors at the same time. To make this algorithm efficient, we identified an effective preconditioner coupled with techniques to generate good initial guesses that significantly accelerates the convergence of the LOBPCG algorithm. We describe in detail our techniques to construct the preconditioner and generate initial guesses, as well as efficient implementation of these techniques on large scale distributed memory clusters.

Although our main focus is to develop, implement and test an efficient LOBPCG approach within MFDn, many of the features of our implementation may be adopted straightforwardly in other configuration interaction (CI) codes. For example, the codes Antoine [10], BIGSTICK [16, 30], KSHELL [31], and NuShellX [8] represent popular CI codes that also solve the large sparse matrix eigenvalue problem arising in nuclear structure calculations. Of these codes, BIGSTICK is also capable of treating three-body interactions and is specifically designed for high-performance computing. The most important difference between BIGSTICK and MFDn is that in MFDn the sparse matrix is stored in core during the diagonalization phase, whereas it is recomputed-on-the-fly in BIGSTICK. As the number of Lanczos iterations, each with a recomputed matrix, can then be time-consuming, BIGSTICK may especially benefit from some of the key features of our LOBPCG approach.

The paper is organized as follows. In Section 2, we review the general formulation of the nuclear many-body problem and the structure of the Hamiltonian matrix to be diagonalized. We describe the LOBPCG algorithm for computing the lowest few eigenvalues of $\hat{H}$ in Section 3 along with a discussion of our recent improvements, which are the main contribution of this paper. In Section 4, we discuss a number of implementation issues. In Section 5, we report the computational performance of the resulting LOBPCG implementation and compare it with the existing Lanczos eigensolver in MFDn.

## 2  The Structure of the Hamiltonian

The efficiency of the block iterative LOBPCG eigensolver hinges, to a large extent, on how efficient the multiplication of $\hat{H}$ with a block of vectors is carried out. Due to its large dimension and the relatively large number of nonzero matrix elements it contains, the matrix $\hat{H}$ is generated in parallel, and stored over a set of distributed memory multiprocessors in a two dimensional fashion, see Section 4.1. The sparsity structures of the local submatrices are not known in advance, and will in general depend on how the many-body basis functions $\Phi_a(\mathbf{r}_1, \mathbf{r}_2, \ldots, \mathbf{r}_A)$ are generated, ordered, and distributed. To achieve good performance, we generate these basis functions in a way to

1. ensure that the parallel sparse matrix vector multiplication is well load-balanced;

2. allow the sparsity structure of $\hat{H}$ to be quickly determined without exhaustively evaluating $\langle\Phi_a|\hat{H}|\Phi_b\rangle$ for every $(a,b)$ pair [22, 35], where $\langle\Phi_a|H|\Phi_b\rangle$ denotes the inner product between $\Phi_a$ and $H\Phi_b$ defined in terms of integration with respect to $\mathbf{r}_1, \mathbf{r}_2, \ldots, \mathbf{r}_A$;

3. enhance data locality by keeping the matrix in block sparse format [1, 2, 9].

All these issues are closely related, and depend on how the many-body basis functions are enumerated, which we will describe in Section 2.1. In Section 2.2, we describe a mechanism for grouping many-body basis states to create matrix tiles (a tile is a rectangular block, i.e., a submatrix, of the matrix $\hat{H}$) that allow us to reduce the amount of work required to probe the sparsity structure of $\hat{H}$ and to enhance data locality when $\hat{H}$ is multiplied with a number of vectors. While the general matrix tiling strategy was presented in [2], our previous approach did not leverage the Hamiltonian structure efficiently during matrix construction; it rather focused on exploring the benefits of using sparse matrix multiple vector multiplication (SpMM) over a series of sparse matrix vector multiplies (SpMV). In this paper, we present a novel grouping strategy to enhance data locality and incorporate this strategy directly during the construction of $\hat{H}$.

## 2.1 Wavefunction Indexing and Enumeration

Each many-body basis function $\Phi_a(\mathbf{r}_1, \mathbf{r}_2, \ldots, \mathbf{r}_A)$ that we use to expand the nuclear wavefunctions and express the Hamiltonian can be indexed by an $A$-tuple of single-particle quantum numbers, i.e., $a = (a_1, a_2, \ldots, a_A)$. We often refer to $\Phi_a(\mathbf{r}_1, \mathbf{r}_2, \ldots, \mathbf{r}_A)$ or simply $a$ as a *many-body basis state*. We will call each component of $a$ that appears in (2), i.e., $\phi_{a_i}(\mathbf{r})$, a *single-particle state*. In nuclear physics we have two types of fermions, protons and neutrons, so that each many-body basis state can actually be written as a product of two Slater determinants of the form (2), one for $Z$ protons and one for $N$ neutrons. Since we are dealing with spin-$\frac{1}{2}$ particles, each single-particle state $\phi_{a_i}(\mathbf{r})$ can either be occupied or unoccupied, but we cannot have more than one particle of the same type in any given single-particle state due to the Pauli exclusion principle.

In MFDn, we use a spherically-symmetric basis, for which the angular part of the single-particle wavefunctions can be expanded in spherical harmonics $Y_l$ (a spherical tensor of rank $l$, with $2l+1$ spherical components). These spherical harmonics are multiplied with radial wavefunctions $R_{nl}(r^2)$ to form the single-particle basis function. Using this expansion, we can label each single-particle basis function $\phi_{a_i}(\mathbf{r})$ by their quantum numbers $(n, l, j, m_j)$: the radial quantum number $n$, counting the number of nodes (excluding the nodes at the origin and the infinity) in the radial wavefunction $R_{nl}(r^2)$; the orbital angular momentum quantum number $l$; the total single-particle angular momentum $j = |l \pm \frac{1}{2}|$ which is the intrinsic nucleon spin $s = \frac{1}{2}$ coupled to $l$; and $m_j$, the magnetic projection of $j$, with $m_j = -j, -j+1, \ldots, j-1, j$. For a spherically-symmetric one-body Hamiltonian, its single-particle states that differ only in the $m_j$ component of the quantum numbers are degenerate. We refer to the set of single-particle states with the same $(n, l, j)$ but different values of $m$ as orbitals; each orbital can be occupied by at most $2j+1$ identical particles.

The nuclear Hamiltonian preserves (among other things) the total angular momentum $J$, the total angular momentum projection $M$, and the parity $P = \pm 1$. MFDn is a so-called $M$-scheme code, in which the many-body basis states (i.e., configuration space) have a specific value of $M$, as well as a specific value of $P$, but not a specific value of $J$. The angular momentum projection is a simple additive quantum number

$$M = \sum_{i=1}^{A} m_j(a_i),$$

and the parity is a multiplicative quantum number

$$P = \prod_{i=1}^{A} p(a_i)$$

with $p(a_i) = (-1)^{l(a_i)}$.[2] In an $M$-scheme basis, the Hamiltonian is block-diagonal with each block characterized by its $M$ and $P$ values. In practice, we perform our calculation for one block at a time, that is, for fixed $M$ and $P$. Using $M = 0$ ($\frac{1}{2}$) for even (odd) $A$, one needs only two runs, one for each parity, to determine the (low-lying) spectrum.

When the single-particle states are chosen to be eigenfunctions of a 3D quantum harmonic oscillators, they can be ordered by their energy levels (or eigenvalues), which are given by $(\frac{3}{2}+N)\hbar\omega$ for some energy scale $\omega$ with $N = 2n + l$ starting with $N = 0$. Note the eigenvalues are actually multiple eigenvalues, because each $N > 0$, often referred to as a shell, corresponds to multiple combinations of $n$ and $l$. It is also important to keep in mind that each $(n, l, j)$ orbital has a multiplicity of $2j + 1$, and thus each harmonic oscillator shell has a multiplicity of $(N + 1)(N + 2)$ for each type of nucleon, i.e., for protons and neutrons separately.

Thus, using the harmonic oscillator basis in combination with the $N_{\max}$ truncation, the set of many-body basis states (configurations) can be given by the following set of parameters

1. the number of protons and neutrons, $Z$ and $N$;

2. the quantum numbers $M$ and parity $P$, discussed above;

3. a many-body basis truncation parameter $N_{\max}$, defined below;

4. the energy scale $\omega$ of the harmonic oscillator.

The truncation parameter $N_{\max}$ sets an upper bound for the sum of $2n(a_i) + l(a_i)$ over all nucleons $i = 1, 2, \ldots, A$, i.e.,

$$\sum_{i=1}^{A} \bigl(2n(a_i) + l(a_i)\bigr) \leq N_0 + N_{\max}, \qquad (3)$$

where $N_0$ is the minimum number of quanta for that nucleus given the Pauli exclusion principle. The inequality (3) limits how large each $n(a_i)$ and $l(a_i)$ can become. Implicitly the limit on $N_{\max}$ also limits the underlying single-particle basis. Furthermore, since parity is positive or negative depending on $l$ being even or odd, even values and odd values of $\sum(2n(a_i) + l(a_i))$ correspond to configurations with opposite parity. Thus, implicitly $N_{\max}$ also selects parity when basis states are retained in decrements of two from $N_{\max}$.

We denote the set of configurations (many-body basis states) $\{a\}$ by $\mathcal{A}$. The size of $\mathcal{A}$ gives the dimension of the (sparse) matrix to be diagonalized; this dimension grows rapidly with increasing number of nucleons and with increasing $N_{\max}$; see Fig. 1. In MFDn, we order single-particle states first by their energy levels (or equivalently, harmonic oscillator shell) as given by $n$ and $l$. States in the same shell are further ordered by their angular momentum $j$, starting from the lowest $j$; and within the same orbital, labeled by $(n, l, j)$, we order the single-particle states according to $m_j = -j, -j + 1, \ldots, j - 1, j$. The many-body basis states are first ordered by the values of $\sum_{i=1}^{A}\bigl(2n(a_i) + l(a_i)\bigr)$, and then enumerated in lexicographical order based on this ordering of the single-particle states.

---

[2]Here $m_j(a_i)$ and $p(a_i)$ are used, respectively, to denote the $m_j$ and the parity $p$ of the single-particle state with the quantum numbers $a_i$, and we will use a similar notation for the quantum numbers $n$, $l$, and $j$ as well.

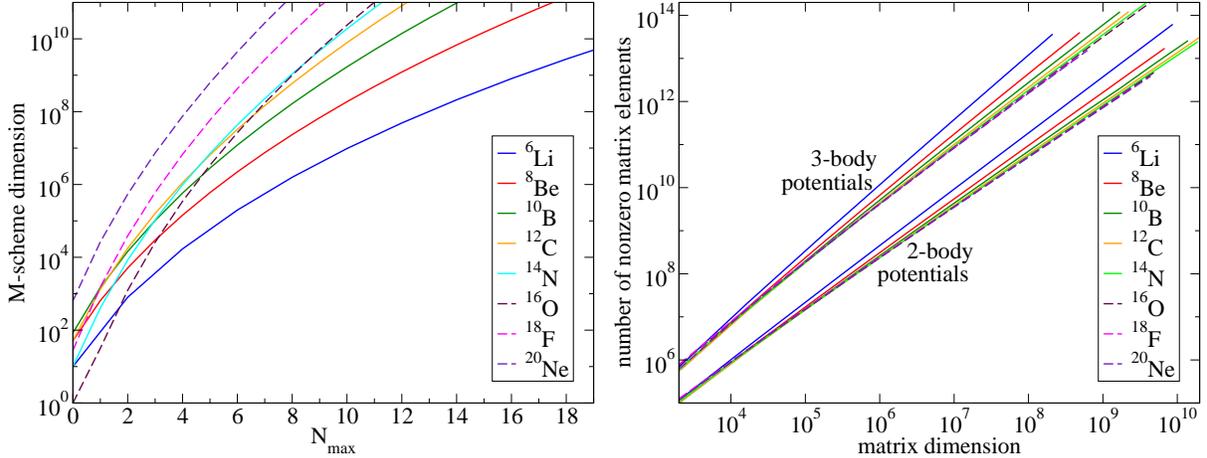

Figure 1: Dimensionality (left) and number of nonzero matrix elements in the lower triangular portion of the symmetric matrix (right) for light nuclei with the same number of protons and neutrons (i.e., $Z = N$), with angular momentum projection $M = 0$ and fixed parity $P$.

## 2.2 Tile and Block Structure of $\hat{H}$

Let us assume that $H$ is an $A$-body Hamiltonian that contains $d$-body interactions, with $d \leq A$. Given two many-body basis states $a$ and $b$ belonging to $\mathcal{A}$, it is evident that

$$\langle \Phi_a | H | \Phi_b \rangle = 0,$$

if $a$ and $b$ differ by more than $d$ single-particle states [11] after summing over the neutron and the proton differences. Therefore, the sparsity of the discretized matrix Hamiltonian $\hat{H}$ can in principle be determined by examining the difference between each pair of $(a, b)$ in $\mathcal{A}$, once all configurations in $\mathcal{A}$ have been generated.

However, going through all pairs of many-body basis states in an exhaustive manner is prohibitively expensive, given the cardinality of $\mathcal{A}$; see Fig. 1. A more efficient scheme used in MFDn to determine the sparsity of $\hat{H}$, is to first consider groups of many-body states based on their orbital quantum numbers $(n, l, j)$ only, without considering $m_j$. That is, many-body basis states that only differ in the $m_j$ values of each single-particle state are placed in the same group labeled by

$$g(a) \equiv \{b \in \mathcal{A} : n(a_i) = n(b_i),\ l(a_i) = l(b_i),\ j(a_i) = j(b_i),\ \text{for } i = 1, 2, \ldots, A\}. \quad (4)$$

This particular choice of grouping is physically natural because the set of many-body states belonging to the same $g(a)$ have the same total single-particle energy when the many-body Hamiltonian contains (or can be approximated by) the one-body Hamiltonian that generates the single-particle states.

The size of these groups $g(a)$'s varies greatly, and can range from just a single $M$-scheme many-body state to many thousands of $M$-scheme many-body states. In order to improve load-balancing, it may be necessary to divide some groups labeled by $g(a)$ into smaller groups, in particular for larger nuclei, with $A > 8$. We can do this by subdividing the orbitals, that is, for example,

$$g(a) \equiv \{b \in \mathcal{A} : n(a_i) = n(b_i),\ l(a_i) = l(b_i),\ j(a_i) = j(b_i),\ m_j^{\min}(a_i) \leq m_j(b_i) \leq m_j^{\max}(a_i)\} \quad (5)$$

and consider these "split orbitals" as the building blocks for the groups of many-body states. Thus, we may control the size of the largest groups.

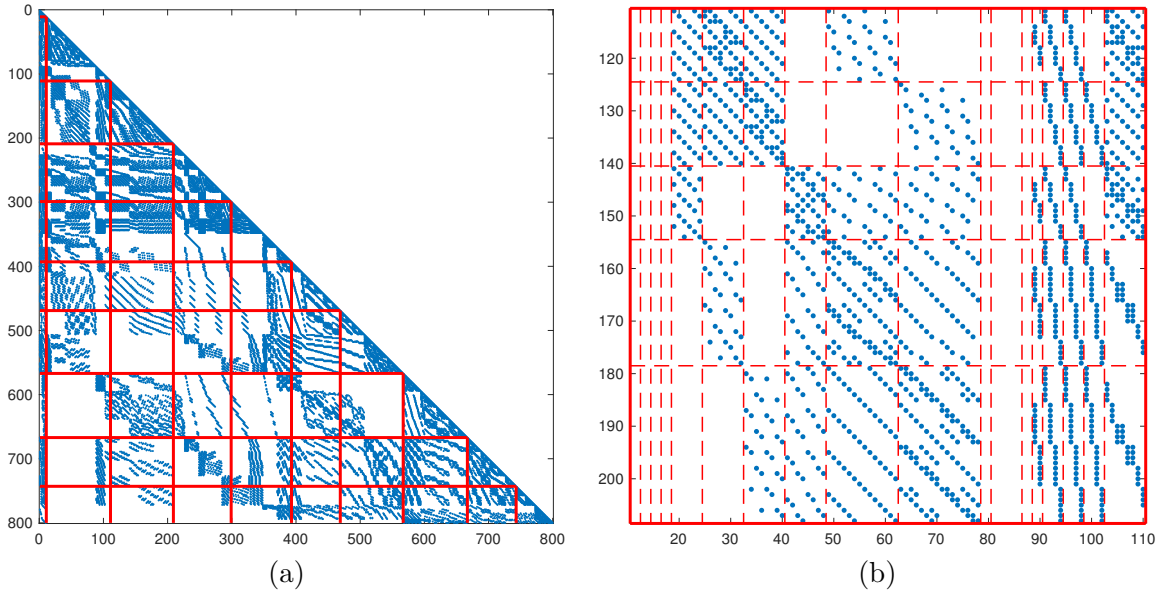

Figure 2: On the left, sparsity structure of $\hat{H}$ for $^6$Li at $N_{\max} = 2$, with 800 many-body basis states, 92 groups of states, 1,826 nonzero tiles and 33,476 nonzero matrix elements in the lower triangular portion of the symmetric matrix. On the right a more detailed plot of one CSB block (the (3,2)-block) of this matrix. Boundaries of CSB blocks and tiles are indicated by the solid and dashed lines, respectively.

The grouping of many-body basis states leads to a partitioning of $\hat{H}$ into many *tiles*. Each tile can be indexed by two group identification labels $g(a) \equiv (g(a_1), g(a_2), \ldots, g(a_A))$ and $g(a') = (g(a'_1), g(a'_2), \ldots, g(a'_A))$, where $g(a_i)$ denotes the $(n, l, j)$ quantum numbers associated with the single-particle state $a_i$. The dimension of the tile is determined by the sizes of $g(a)$ and $g(a')$. In MFDn, we first perform pairwise comparisons of $g(a)$ and $g(a')$ to determine (potentially) nonzero tiles. If $g(a_i) \neq g(a'_j)$ for more than $d$ single-particle states $i$ and $j$, the entire tile indexed by $g(a)$ and $g(a')$ is zero. If $g(a_i) \neq g(a'_j)$ for up to $d$ single-particle states, then the tile generally contains nonzero matrix elements.

Note that the sparsity of these nonzero tiles depends on the number of single-particle states that differ between $g(a)$ and $g(a')$. The *diagonal* tiles, with $g(a) = g(a')$, are generally the least sparse (most dense) tiles, whereas tiles for which $g(a)$ and $g(a')$ differ by exactly $d$ particles are the most sparse (least dense). In principle, additional blocking techniques such as those described in [35] can be used to identify potentially nonzero tiles, or nonzero matrix elements within the tiles, and further reduce the number of pairwise comparisons of (groups of) many-body basis states.

In Fig. 2, we show the sparsity pattern for a small case, $^6$Li, with three protons and three neutrons, truncated to $N_{\max} = 2$. The $M = 0$ basis space dimension is 800, and the number of groups is 92. Thus, there are in principle 4,186 tiles in half of the symmetric matrix, and 320,400 matrix elements. The number of (potentially) nonzero tiles is 1,826, corresponding to a "tile-sparsity" of 0.44, and the total number of nonzero matrix elements within these nonzero tiles is 33,476, corresponding to a sparsity of 0.104.

A well known technique for improving the performance of sparse matrix computation on modern processors is blocking, i.e., partitioning the matrix into a number of smaller blocks that can be loaded into fast memory (e.g., level-2 cache) and processed one block at a time. It is tempting to use the tile structure of $\hat{H}$ to naturally block the matrix. However, this is not the best strategy.

From Fig. 2(b), it is evident that the tiles are not uniform in size. For the small problem size demonstrated in that figure, some tiles are as small as one by one (i.e., containing a single matrix element) whereas the largest one is 52 by 52. The average group size is 8.7 for this particular example. For the problems of interest, which are much larger in basis dimension, the smallest tiles will also be as small as one by one, but the largest tiles can have dimensions of thousands or even more. Keeping track of blocks of small sizes introduces significant overhead and the large variation in block sizes can make load balancing highly nontrivial. Therefore, we introduce an additional level of blocking to group $g(a)$'s into "sets of groups" containing up to $\beta$ many-body basis states for some modest $\beta$ value (e.g., $\beta = 4000$). We shall refer to the resulting blocks as *CSB blocks*, named after the sparse matrix format we adopt as described in Section 4.2. Such a block structure is indicated by the solid vertical and horizontal lines in Fig. 2(a) for $\beta = 100$. Although it is not obvious from the figure, all block boundaries do coincide with tile boundaries. This is ensured by not enforcing a constant CSB block size, and allowing the aggregate size of the sets of groups in each block to be slightly less than $\beta$. In addition, we impose block boundaries between subspaces corresponding to different $N_{\max}$ truncations. In the particular example shown in Fig. 2, $\beta$ is set to 100, and the average size of obtained CSB blocks is 80 by 80. There is also a block boundary between the $N_{\max} = 0$ and $N_{\max} = 2$ subspaces (in this case the $N_{\max} = 0$ subspace has a dimension of 10, the left-most block column).

The ideal value of $\beta$ and corresponding ideal largest group size strongly depend on both the nuclear physics problem (number of nucleons and the specific basis space truncation) and on the computer hardware (available cache sizes, memory bandwidth, number of OpenMP threads etc.), and has been discussed in detail in our previous work [1, 2]. The block partitioned matrix is stored in a compressed sparse block (CSB) format proposed in [9]. We will discuss the CSB data structure in Section 4.

## 3 The LOBPCG Algorithm

Most CI codes used in nuclear physics, including MFDn, use a Lanczos based eigensolver with full orthogonalization to ensure numerical stability. There are three main advantages of using the LOBPCG algorithm over the existing Lanczos solver. First, the LOBPCG algorithm allows us to make use of preconditioners to be discussed below. With good preconditioners, the convergence of the algorithm can be accelerated significantly as we demonstrate in Section 5. Second, the LOBPCG algorithm can make use of good approximations to different eigenvectors as the starting guess. In contrast, the Lanczos algorithm can only use one of these vectors as the starting vector. Finally, the implementation of the LOBPCG algorithm naturally depends on kernels with high arithmetic intensities, i.e., multiplication of a sparse matrix with multiple vectors (SpMM) and level-3 BLAS operations to update the block vectors.

In this section, we first review the basic LOBPCG algorithm and then discuss how to make use of the special properties of the Hamiltonian in MFDn to accelerate the eigensolver.

### 3.1 The basic LOBPCG algorithm

We denote the eigenvalues of an $D \times D$ discretized Hamiltonian $\hat{H}$ by $\lambda_1 \leq \lambda_2 \leq \cdots \leq \lambda_D$ and the corresponding eigenvectors by $x_1, x_2, \ldots, x_D$. It is well known that $X_* = [x_1, x_2, \ldots, x_k]$ is the solution to the trace minimization problem

$$\min_{X^T X = I_k} \operatorname{trace}(X^T \hat{H} X). \tag{6}$$

The locally optimal block preconditioned conjugate gradient (LOBPCG) algorithm developed by Knyazev [17] seeks the solution of (6) by updating the eigenvector approximation $X^{(i)}$ according to

$$X^{(i+1)} = X^{(i)}G_1 + R^{(i)}G_2 + P^{(i)}G_3, \qquad (7)$$

where $R^{(i)}$ is the projected gradient

$$R^{(i)} = \left[I - X^{(i)}(X^{(i)})^T\right]\hat{H}X^{(i)}, \qquad (8)$$

and $P^{(i)}$ is the previous search direction defined recursively as

$$P^{(i)} = R^{(i-1)}G_2 + P^{(i-1)}G_3,$$

with $P^{(0)} = 0$. The superscript $(i)$ denotes the iteration number. The $k \times k$ matrices $G_1$, $G_2$, and $G_3$ are chosen in each iteration $i$ to minimize the objective function in (6) within the subspace spanned by columns of

$$Y \equiv \begin{bmatrix} X^{(i)} & R^{(i)} & P^{(i)} \end{bmatrix}.$$

The solution to this subspace minimization problem can be obtained from the solution of the following generalized eigenvalue problem

$$(Y^T \hat{H} Y)G = (Y^T Y)G\Theta, \qquad (9)$$

where

$$G = \begin{bmatrix} G_1 \\ G_2 \\ G_3 \end{bmatrix}$$

consists of eigenvectors associated with the smallest eigenvalues of the matrix pencil $(Y^T \hat{H} Y, Y^T Y)$. These eigenvalues appear in the diagonal matrix $\Theta$.

To accelerate the convergence of the algorithm, preconditioning is adopted in practice. More precisely, let $W^{(i)} = T^{-1}R^{(i)}$, where $T$ is the *preconditioner*. Then

$$Y \equiv \begin{bmatrix} X^{(i)} & W^{(i)} & P^{(i)} \end{bmatrix}$$

is used to construct the projection subspace. We also need to correspondingly replace $R^{(i)}$ by $W^{(i)}$ in (7). The main steps of the LOBPCG algorithm are outlined in Algorithm 1. A practical implementation of the LOBPCG algorithm requires a more careful treatment of the subspace span($Y$). For instance, it is advisable to choose $k$ somewhat larger than the actual number of desired eigenpairs. Orthogonalization on $Y$ is also important for numerical stability. We refer the readers to [12, 14] for techniques on a robust implementation.

In the rest of this section, we discuss how to make use of the special properties of the nuclear Hamiltonian in the context of MFDn to accelerate the eigensolver.

## 3.2 Preconditioner

The convergence rate of the LOBPCG algorithm is related to the condition number of $\hat{H}$. The method can be accelerated when a good preconditioner is chosen in step 6 of Algorithm 1. From an optimization point of view, the use of a symmetric positive definite preconditioner $T = LL^T$ where $L$ is a lower triangular Cholesky factor of $T$, can be considered as a change of variable $Y = L^{-1}X$ applied to (6) to yield a preconditioned minimization problem

$$\min_{Y^T(L^T L)Y = I_k} \text{trace}(Y^T L^{-T} \hat{H} L^{-1} Y). \qquad (10)$$

**Algorithm 1** The locally optimal block preconditioned conjugate gradient (LOBPCG) algorithm.

**Input:** $\hat{H} = \hat{H}^T \in \mathbb{R}^{D \times D}$, $X^{(0)} \in \mathbb{R}^{D \times k}$ satisfying $(X^{(0)})^T X^{(0)} = I_k$.
**Output:** The $k$ smallest eigenpairs of $\hat{H}$, i.e., $X \in \mathbb{R}^{D \times k}$ and $\Theta \in \mathbb{R}^{k \times k}$ such that $\hat{H}X = X\Theta$.
1: Set $\Theta^{(0)}$ to the Rayleigh quotient of $X^{(0)}$.
2: Set $P^{(0)} = []$.
3: **for** $i = 1, 2, \ldots$ **do**
4:   Set $R^{(i)} = \hat{H}X^{(i)} - X^{(i)}\Theta^{(i)}$.
5:   Check convergence.
6:   Set $W^{(i)} = T^{-1}R^{(i)}$.
7:   Orthogonalize $\{X^{(i)}, W^{(i)}, P^{(i)}\}$ and solve the projected eigenvalue problem on the subspace span $\{X^{(i)}, W^{(i)}, P^{(i)}\}$.
8:   Set $\Theta^{(i+1)}$ to the $k \times k$ diagonal matrix with the $k$ smallest Ritz values as its diagonals.
9:   Choose the Ritz vectors correspond to the $k$ smallest Ritz values as $X^{(i+1)}$.
10:  Orthogonalize $X^{(i+1)} - X^{(i)}$ against $X^{(i)}$ as $P^{(i)}$.
11: **end for**

If the condition number of $L^{-T}\hat{H}L^{-1}$ is smaller than that of $\hat{H}$ while $T = LL^T$ itself is not too ill-conditioned, then the search directions generated in LOBPCG are likely to point more directly towards the global minimizer of (10) [13].

Because the LOBPCG algorithm is closely related to the inexact Newton type methods such as the Jacobi–Davidson algorithm [32], the use of a preconditioner can be viewed as a way to include inexact Newton correction in a search space from which approximate eigenpairs of $\hat{H}$ are extracted. Such a viewpoint allows the choice of $T$ to be made more flexibly. In particular, we can choose different $T$'s for different columns of $R$ in step 6 of Algorithm 1. We use $T_j$ to denote these different choices of $T$. Each $T_j$ can be chosen to approximate $\hat{H} - \mu_j I$, where $\mu_j$ is an approximation to $j$th eigenvalue $\lambda_j$. Theoretical analysis for this type of shift-and-invert based preconditioners can be found in [5, 25, 26].

In practical computations the preconditioning strategy is often problem dependent. When choosing the preconditioner $T$, we need to balance the objectives of convergence acceleration with that of keeping the cost of step 6 of Algorithm 1 low. We would like to choose a preconditioner that can significantly reduce the number of LOBPCG iterations required to reach convergence without introducing significant additional computational cost per iteration. The types of preconditioners we have examined include the following:

- $T_j = K - \mu_j I$, where $K$ is the kinetic energy part of $\hat{H}$. To apply such a preconditioner, we need to solve linear equations of the form $(K - \mu_j I)z_j = r_j$, for $j = 1, \ldots, k$. Because $K - \mu_j I$ cannot be efficiently factored, these linear equations need to be solved iteratively. Each step of the iterative solution requires multiplying $K$ with multiple vectors. Because the cost of such a multiplication is not significantly lower than that of multiplying $\hat{H}$ with multiple vectors, the cost of applying the preconditioner is relatively high. As a result, the use of such a preconditioner typically does not result in an overall reduction in computational time although the number of LOBPCG iterations required to reach convergence can be reduced slightly. It is possible to drop small entries of $K$ to make the preconditioner sparser, thereby reducing the cost of applying the preconditioner. However, it is not easy to determine an appropriate threshold for dropping entries of $K$. We found that dropping too many entries of $K$ can actually result in an increase in the number of LOBPCG iterations.

- $T_j = \text{diag}\left\{\hat{H}_{N_{\max}-2} - \mu_j I, I\right\}$ where $\hat{H}_{N_{\max}-2}$ is the Hamiltonian constructed from a smaller configuration space. Although this type of preconditioner makes intuitive sense, our numerical experiments indicate that it is not effective, and often results in an increase in the number of LOBPCG iterations.

- $T_j = \tilde{H} - \mu_j I$, where $\tilde{H}$ is a block diagonal approximation of $\hat{H}$. Because the diagonal blocks of $\hat{H}$ contains matrix elements that correspond to Slater determinants that interact strongly through the many-body Hamiltonian, we expect $\tilde{H}$ to capture main spectral properties of $\hat{H}$. When the sizes of the diagonal blocks are reasonably small, this preconditioner is cheap to apply. Several block partitioning strategies can be used to create the diagonal preconditioner. For example, we can use the diagonal tiles, diagonal CSB blocks, or the local Hamiltonian assigned to each diagonal processor (see Section 4.1) to construct this preconditioner. In the extreme case, we can also use just the diagonal elements of $\hat{H}$ to construct a diagonal preconditioner. If we use $\tilde{H} = \text{diag}\left\{\tilde{H}_{[1]}, \ldots, \tilde{H}_{[m]}\right\}$ to denote the block diagonal preconditioner, where $m$ is the number of diagonal blocks and varies based on the particular strategy used, then applying the preconditioner amounts to solving $m$ independent linear systems of equations of the form

$$(\tilde{H}_{[i]} - \mu I)W_{[i]} = R_{[i]} \tag{11}$$

for $i = 1, \ldots, m$, where $W_{[i]}$ and $R_{[i]}$ are, respectively, the $i$th chunk of $W$ and $R$ partitioned conformally with the partitions $\tilde{H}_{[i]}$. When the block size of $\tilde{H}_{[i]}$ is very small, (11) can be solved by a direct method. Otherwise, we can use an iterative method such as the *full orthogonal method* (FOM) [27] or the *generalized minimal residual method* (GMRES) [29]. Our experiments show that two or three iterations are typically sufficient to accelerate the convergence of LOBPCG. Using larger diagonal $\tilde{H}_{[i]}$ blocks may results in slightly fewer LOBPCG iterations, but this clearly comes at the expense of increased computational costs during the FOM or GMRES iterations which require multiplications of $\tilde{H}_{[i]}$ with multiple vectors. We found that a good compromise between reducing the number of LOBPCG iterations and keeping the cost of applying the preconditioner low per iteration is to use the diagonal tiles associated with many-body groups in the block diagonal approximation. To achieve faster convergence for LOBPCG, we allow different shifts $\mu_j$ for different right hand sides in (11).

## 3.3 The initial guess for LOBPCG

In principle, almost any linearly independent set of vectors can be used as the initial guess of LOBPCG. However, when good approximations to several eigenvectors are available, the LOBPCG algorithm allows us to use all these approximate eigenvectors as the initial guess. Therefore it is beneficial to find good approximate eigenvectors before applying the eigensolver. Fortunately, in MFDn, there are a number of efficient ways to construct reasonably good approximate eigenvectors. In the following we discuss several choices, all based on approximations to the Hamiltonian. Some techniques can be used not only for MFDn, but also for more general configuration interaction calculations.

### 3.3.1 Constructing from a smaller configuration space

Let $\hat{H}$ be partitioned as

$$\hat{H} = \begin{bmatrix} \hat{H}_{11} & \hat{H}_{12} \\ \hat{H}_{21} & \hat{H}_{22} \end{bmatrix},$$

where $\hat{H}_{11}$ is the discretized Hamiltonian corresponding to a smaller configuration space obtained by $\sum_{i=1}^{A}(2n(a_i) + l(a_i)) \leq N_0 + N_{\max} - 2$. As the dimension of $\hat{H}_{11}$ is much smaller compared to that of $\hat{H}$, computing the $k$ smallest eigenpairs of $\hat{H}_{11}$ requires relatively low cost. Suppose that we have already computed $\hat{H}_{11}X_1 = X_1\Theta_1$, for instance, by using the LOBPCG algorithm described above. Approximate eigenvectors of $\hat{H}$ can then be obtained by padding $X_1$ with zeros, i.e.,

$$X^{(0)} = \begin{bmatrix} X_1 \\ 0 \end{bmatrix}. \tag{12}$$

Then

$$\frac{\|r_j^{(0)}\|}{\|\hat{H}\|\|x_j^{(0)}\|} = \frac{\|\hat{H}x_j^{(0)} - x_j^{(0)}\theta_j\|}{\|\hat{H}\|\|x_j^{(0)}\|} = \frac{\|\hat{H}_{21}x_j^{(0)}\|}{\|\hat{H}\|\|x_j^{(0)}\|} \leq \frac{\|\hat{H}_{21}\|}{\|\hat{H}\|} =: \epsilon.^3$$

The ratio $\epsilon$ is much less than one when the truncation error of the smaller configuration space is small. Thus (12) provides a good initial guess for LOBPCG.

An alternative but computationally more expensive approach is to set

$$X^{(0)} = \begin{bmatrix} X_1 \\ Y_1 \end{bmatrix}, \tag{13}$$

where the $j$th column of $Y_1$ is chosen as the $j$th column of $(\hat{H}_{22} - \theta_j I)^{-1} X_1$. The relative residual norm is then bounded by

$$\frac{\|r_j^{(0)}\|}{\|\hat{H}\|\|x_j^{(0)}\|} = \frac{\|\hat{H}_{21}^T(\hat{H}_{22} - \theta_j I)^{-1}\hat{H}_{21}x_j^{(0)}\|}{\|\hat{H}\|\|x_j^{(0)}\|} \leq \frac{\|\hat{H}\|}{\mathrm{gap}_j}\epsilon^2 = O(\epsilon^2),$$

where $\mathrm{gap}_j := \|(\hat{H}_{22} - \theta_j I)^{-1}\|^{-1}$. Therefore, (13) is in general better than (12). The price to pay is that (13) requires solving linear systems to obtain $Y_1$.

In practice, we adopt (12) as the initial guess for its lower computational cost. However, we remark that the block diagonal preconditioning provides us an opportunity to achieve better approximation quality. Suppose that $T_j = \mathrm{diag}\left\{\hat{H}_{11}, \hat{H}_{22}\right\} - \theta_j I$ is used as the preconditioner for $r_j^{(0)}$ with the initial guess (12). Then it can be verified that

$$W^{(0)} = -\begin{bmatrix} 0 \\ Y_1 \end{bmatrix}.$$

As the alternative initial guess (13) belongs to span $\{X^{(0)}, W^{(0)}\}$, one step of LOBPCG achieves the $O(\epsilon^2)$ accuracy. Although, compared to $T_j$, the tile-partitioned preconditioner leaves out more off-diagonal blocks from $\hat{H} - \theta_j I$, we still expect some accuracy improvement of the approximate eigenvectors after one step of LOBPCG.

Finally, we remark that when performing LOBPCG on a smaller configuration subspace, we can actually utilize LOBPCG recursively on an even smaller configuration space to produce an initial guess. Such a recursive use leads to a *multi-level* LOBPCG algorithm for CI calculations. In MFDn, we choose to perform three LOBPCG runs, corresponding to configuration spaces specified by $N_{\max} - 4$, $N_{\max} - 2$, and $N_{\max}$. For the smallest configuration space specified by $N_{\max} - 4$ we find that LOBPCG with random initial guess is adequate for our purpose.

---

[3] We use the lower case variables, e.g., $x_j$, to denote the $j$th column of $X$.

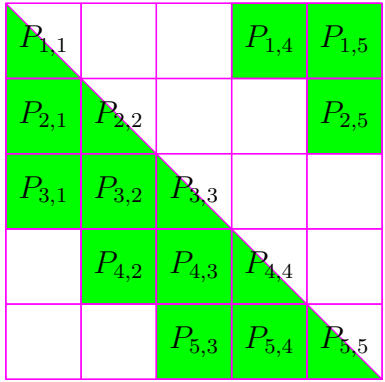

Figure 3: A triangular processor grid with $n_d = 5$ using the BCM mapping in the current implementation of MFDn. Highlighted matrix blocks are stored on corresponding processors.

### 3.3.2 Use the eigenvectors of a closely related Hamiltonian

In addition to the multi-level LOBPCG method discussed above, MFDn also allows other types of initial guesses to be used for LOBPCG. These types of initial guesses include:

- eigenvectors of a simplified Hamiltonian that includes an external field;
- eigenvectors of a Hamiltonian with a different type of nuclear potential or a different set of parameters.

We shall discuss them in more detail in Section 5.

## 4 A High Performance LOBPCG Implementation

In this section, we briefly describe a number of implementation issues that are critical for achieving high performance. In Section 4.1 we summarize the partitioning of the Hamiltonian and LOBPCG vectors, as well as the collective communication operations involved in a distributed memory setting. The algorithms used for this purpose are similar to our Lanczos based implementation which is described in more detail in our previous work [3, 4]. Then in Sections 4.2 through 4.4, we describe details specific to our LOBPCG implementation.

### 4.1 The distribution of $\hat{H}$ and vector segments

Since $\hat{H}$ is symmetric and the number of nonzero elements in $\hat{H}$ increases rapidly as the dimension of $\hat{H}$ increases, we store only the lower triangular part of the Hamiltonian. We would like to distribute $\hat{H}$ among different processing units in such a way that each processing unit will perform roughly the same number of operations in parallel sparse matrix computations. In [3] we developed a grouping and mapping scheme, which we refer to as *balanced column major* (BCM) mapping. Fig. 3 shows the placement of 15 processors on a $5 \times 5$ processor grid using the BCM mapping. The BCM scheme limits the number of processors in each row or column communication group to $(n_d + 1)/2$.[4] We refer to [3] for detailed discussions.

The distribution of dense blocks of LOBPCG vectors requires two different schemes, depending on the type of operations applied to these vectors. When performing the SpMM operation $U = \hat{H}W$, we partition the dense blocks of vectors $W$ and $U$ conformal to the partition of $\hat{H}$, and distribute them among diagonal processors. Under the BCM scheme, we denote the submatrix of $\hat{H}$ assigned to the processor $P_{i,j}$ by $\hat{H}_{i,j}$. Since we store only half of the symmetric Hamiltonian, each processor is responsible for sparse matrix computations related to both the submatrix assigned to it, as well as the transpose of this submatrix. Then each processor needs to performs two local SpMM's of the form $U_i = \hat{H}_{i,j} W_j$ and $U_j = \hat{H}_{i,j}^T W_i$. Once local SpMM operations are completed, two reductions are required to aggregate partial products $U_i$ and $U_j$ into the global result vector $U$.

---

[4]Here we require $n_d$ to be an odd number so that $(n_d + 1)/2$ is an integer.

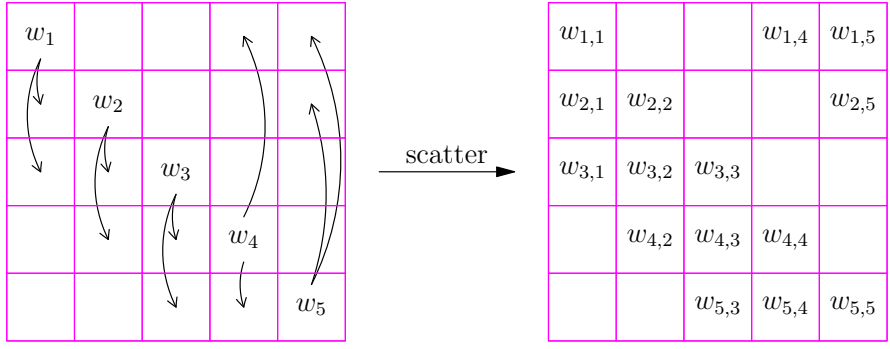

Figure 4: Scatter a distributed vector $w$ from $n_d$ diagonal processors ($n_d = 5$ here) to all processors. Each diagonal processor scatters its portion of $w$ to $(n_d + 1)/2$ processors belonging to the same column processor group.

Although the multiplication of $\hat{H}$ with a block of vectors constitutes the major computational work in each LOBPCG iteration, the cost of other linear algebra operations, such as the projected gradient computation defined in (8), cannot be ignored in a high performance implementation. If $W$ and $U = \hat{H}W$ are partitioned by rows and distributed among $n_d$ diagonal processors only, the number of processors that can effectively be used to perform the dense matrix–matrix multiplications required in (8) is limited to $n_d$. To maximize the level of concurrency, we further divide each vector block of $W_i$ and $U_i$ on $P_{i,i}$ and scatter these vector segments among $(n_d + 1)/2$ processors in the column communication group of $P_{i,i}$; see Fig. 4. As a result, the dense matrix–matrix multiplications of (7), (8), and (9) can be carried out using all $n_d(n_d + 1)/2$ processors available.

To summarize, the data distribution scheme discussed above means that each LOBPCG step involves several distinct collective communication patterns, as illustrated in Fig. 5, as well as local computation. Schematically, we have

1. Before performing $U = \hat{H}W$, the segments of $W$ which are distributed among column processor groups are gathered onto the diagonal processor within each column communication group;

2. The diagonal processors broadcast the gathered submatrices across their row and column communication groups in preparation for the distributed SpMM computations. As a result, each processor (with the exception of the diagonal processors) receives two submatrices of $W$.

3. After processor $P_{i,j}$ performs two SpMM's, one product is reduced among its column communication group onto the diagonal processor $P_{j,j}$, and the other one is reduced among the row communication group onto the diagonal processor $P_{i,i}$.

4. Each diagonal scatters the final result among $(n_d+1)/2$ processors within its column communication group in preparation for the parallel dense matrix–matrix multiplications of LOBPCG shown in Eqs. (7), (8), and (9).

We remark that in practice we can perform step 1 and the broadcast among the column communication group of step 2 by a single collective MPI call. Next, we can overlap the communication across the row communication group in steps 2 and 3 with the local SpMM and transpose SpMM operations. And finally we can perform the reduction operation among the column communication group of step 3 and the subsequent scatter by a single collective MPI call.

## 4.2 Data structures for storing submatrices of $\hat{H}$ and local vector blocks

As described above, the distribution of the $\hat{H}$ matrix, block vectors and collective communication patterns in our LOBPCG implementation are similar to those of the Lanczos solver from the earlier versions of MFDn. A major difference between the two eigensolvers is that LOBPCG requires working with multiple vectors, as opposed to the single vector iterations performed in the Lanczos algorithm. In fact, this difference leads to the third advantage of the LOBPCG over Lanczos discussed in Section 3.1—with LOBPCG, significant performance improvements can potentially be achieved by exploiting the high arithmetic intensity of the main computational kernels, namely the sparse matrix computations and the update of LOBPCG vector blocks. As we shall show in Section 5, despite the use of an effective preconditioner and good initial guesses, the total number of sparse matrix vector operations required by LOBPCG can still be significantly higher than that of Lanczos. Therefore the choice of data structures that will allow us to effectively exploit this third advantage is critical for the overall performance of LOBPCG.

To store the local block of vectors in LOBPCG, we adopt the row major storage. Compared to the column major alternative, row major storage is more convenient for the collective communications described in the previous subsection because in this scheme, data to be transmitted to other processes are stored contiguously in memory. More importantly, row major storage of vector blocks provides good data locality for performing SpMM—for each nonzero matrix element, several required entries of the vector blocks can be loaded into the cache with a single load. This is especially important for modern architectures with limited cache space per core, and a widening performance gap between the processor and the memory system.

Compressed sparse row (CSR) format [33] is one of the most commonly used data structures for sparse matrices, largely due to its simplicity and good performance in general. However, for MFDn, it has two important pitfalls. First, the CSR format cannot exploit the tiled sparsity structure in local MFDn matrices, which leads to poor cache utilization. The second and more significant issue is that performing the multiplication $\hat{H}_{i,j}^T W_i$ in parallel within a single MPI task requires accumulating the elementwise products in an array that must be shared and synchronized among different threads. The synchronization overhead severely degrades the performance of this multiplication [2].

To achieve high performance while retaining a small memory footprint, we use the compressed sparse block (CSB) format proposed by Buluç *et al.* [9] to store $\hat{H}_{i,j}$ submatrices locally. Each CSB block is independently addressable through a 2D array of pointers. The block partitioning in CSB is an excellent match with the tiled structure of MFDn matrices that results from the grouping of many-body states. The CSB format does not specify how nonzeros within each block are stored—we use the sparse coordinate format where the local row and column offsets for each nonzero element are stored along with the nonzero value. An important requirement to ensure a small memory footprint in this case is to keep $\beta < 2^{16}$, so that 16-bit (instead of 32-bit) indexing can be used for row and column offsets. For MFDn matrices, we empirically determined the ideal value of $\beta$ to be from 2,000 to 4,000 depending on the size of the vector blocks in LOBPCG.

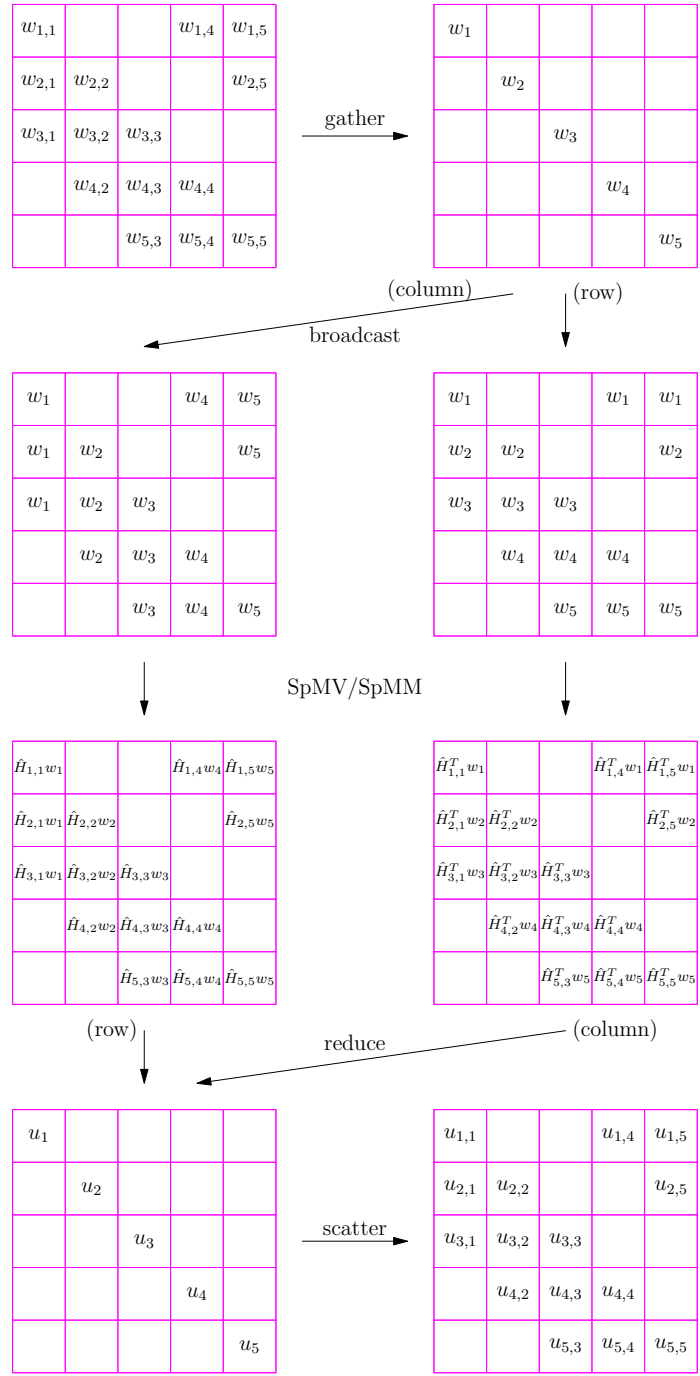

Figure 5: Collective communication in SpMV/SpMM operations.

## 4.3 The construction and application of the preconditioner

As discussed in Section 3.2, we prefer tile-partitioned block diagonal preconditioners in MFDn in which the diagonal tiles of the $\hat{H}$ matrix are used. Hence, when constructing $\hat{H}$, we set up a separate sparse matrix $\tilde{H}$ that consists of these diagonal tiles, and distribute it across the processor grid conformal to the distribution of the vector blocks described in Section 4.1. Consequently, all linear systems in the preconditioning stage are completely local, and can be solved without MPI communication.

Within each MPI process, the local linear system can be further decoupled according to the block diagonal matrix structure of the preconditioner. These decoupled linear systems are assigned to different threads and solved in parallel in our hybrid MPI/OpenMP parallel implementation. Furthermore, since each linear system contains multiple right hand sides, multiple calls to the FOM solver are needed to obtain the approximate solution. To increase arithmetic intensity, we fuse these FOM calls together so that in each iteration we can use a single SpMM instead of a number of SpMV's to construct different Krylov subspaces associated with different right-hand sides at the same time. Up to three FOM iterations are performed for all right-hand sides.

Finally, we remark that we make several conservative choices for the preconditioning to enhance the robustness:

1. No preconditioning is applied for the first three LOBPCG iterates in any case;

2. We do not apply the preconditioners if the relative residual norm of $(\lambda_1, x_1)$ is above $10^{-1}$;

3. If $(\lambda_i, x_i)$ is nearly converged, we set $\mu_i = \theta_i - 2\|r_i\|/\|x_i\|$;

4. If $(\lambda_i, x_i)$ is far from converged while $(\lambda_1, x_1)$, …, $(\lambda_{i-1}, x_{i-1})$ are not, we set $\mu_i = \cdots = \mu_k = \mu_{i-1}$.

## 4.4 A hybrid MPI/OpenMP parallel implementation

There are mainly four types of matrix operations in the LOBPCG algorithm:[5]

1. Preconditioning $W \leftarrow \hat{H} R$;

2. Sparse matrix–matrix multiplication $U \leftarrow \hat{H} W$;

3. Block inner product $X^T Y$;

4. Linear combination $X \leftarrow YG$.

In addition to the concurrency introduced by distributing the ingredients of the eigensolver among different processors, thread level parallelism available through OpenMP can be added to all computations performed within each MPI process. Exploiting shared memory parallelism reduces the internode communication volume [4], as well as the total memory footprint of MFDn [34].

---

[5]As a remark, in MFDn we perform the first two types of operations in single precision because $\hat{H}$ is stored using single precision floating point numbers, while the latter two types of operations are performed in double precision to enhance the numerical stability of LOBPCG.

**Preconditioning** In the preconditioning stage, the computational work involves solving a set of small linear systems approximately. As these small linear systems are completely independent, they can be assigned to different threads and solved in parallel. In our implementation, we exploit this thread level parallelism, and let each thread perform a sequential (approximate) linear solve with multiple right-hand sides. We adopt OpenMP's dynamic scheduling option to resolve potential load imbalances due to different sizes of the linear systems involved.

**Sparse Matrix Multiple Vector Multiplication (SpMM)** The SpMM operation requires performing both $\hat{H}_{i,j}W_j$ and $\hat{H}_{i,j}^T W_i$. As mentioned above, each CSB block that stores nonzero matrix elements is independently addressed through a 2D array of pointers. Hence, the multiplication $\hat{H}_{i,j}W_j$ can be performed by processing this 2D array by rows, while the multiplication $\hat{H}_{i,j}^T W_i$ can be performed by processing the same 2D array by columns. The CSB format automatically enables cache blocking for each tile, thereby allowing us to achieve excellent performance [2].

For thread parallelism, we parallelize the loops over the 2D pointer array to perform $\hat{H}_{i,j}W_j$ and $\hat{H}_{i,j}^T W_i$ with OpenMP. Let us assume that the local matrix $\hat{H}_{i,j}$ has been partitioned into $\ell_i$ row blocks, corresponding to $\ell_i$ groups. Each row block contains $\ell_j$ column blocks, also conformal to the tile partition along the columns. When performing $\hat{H}_{i,j}W_j$, different threads handle different row blocks in parallel. As the size of row blocks and the number of nonzeros contained within each block vary, we adopt dynamic scheduling to improve load balance among different threads. Similarly, when performing the transpose operation $\hat{H}_{i,j}^T W_i$, thread level parallelism is exploited at the level of column blocks. The ability to easily traverse the CSB blocks in column major order allows us to avoid any race conditions while performing the transpose SpMM operation. Unlike the case with CSR storage, we achieve almost identical performance for $\hat{H}_{i,j}W_j$ and $\hat{H}_{i,j}^T W_i$ with a thread-parallel CSB implementation [2].

**Block inner product and linear combination** The dense matrix–matrix operations of LOBPCG given in Eqs. (7), (8), and (9) that involve vector blocks $W^{(i)}$, $R^{(i)}$, and $P^{(i)}$ can be described as block inner products and linear combinations. For the purposes of this subsection, we refer to these vector blocks as $X$ and $Y$, generically.

The block inner products $X^T Y$ are parallelized across MPI processes as local block inner products $X_i^T Y_i$ followed by a reduce operation. The local block inner product operation involves two (extremely) tall-and-skinny matrices $X_i$ and $Y_i$. While high performance BLAS libraries can be used for this purpose, Intel's MKL and Cray's LibSci libraries exhibit poor performance, most likely due to the unusual shape of the matrices involved. Therefore we implemented a custom thread-parallel kernel for block inner products, where we further partition the local tall-and-skinny matrices into row blocks, and make several sequential library GEMM calls with these small rectangular matrices. Each small sequential GEMM call is assigned to a different thread using dynamic scheduling, and partial results from each thread is aggregated at the end, again using thread-parallelism. By doing so, we achieve 1.7× to 10× speedup over performing the block inner products as direct calls to the library GEMM implementations [1].

The operation of linear combination $Y = XG$, where $G$ is a small square matrix whose rank is equal to the number of columns in $X$, is simpler than the block inner product because no interprocess MPI communication is needed. Similar to the local block inner product, we implemented a custom thread-parallel linear combination kernel that partitions the local tall-and-skinny matrix $X_i$ into row blocks and assigns the small linear combination operations to different threads. We observed only small gains in this case over the default GEMM implementation in MKL and LibSci [1].

Table 1: List of the test cases. Top four cases are chosen with the two-body interaction; bottom four cases include a three-body interaction.

| Nuclei | $N_{\max}$ | interaction | dimension of $\hat{H}$ | # of nonzeros in $\hat{H}$ | # of cores |
|---|---|---|---|---|---|
| $^9$Li | 8 | 2-body | 37,145,993 | 24,354,119,738 | $28 \times 6$ |
| $^{10}$B | 8 | 2-body | 159,953,284 | 123,638,482,768 | $120 \times 6$ |
| $^8$Be | 10 | 2-body | 187,304,858 | 206,509,001,984 | $276 \times 6$ |
| $^7$Li | 12 | 2-body | 252,281,462 | 400,190,275,118 | $496 \times 6$ |
| $^6$Li | 8 | 3-body | 1,578,624 | 23,402,203,878 | $28 \times 6$ |
| $^8$He | 8 | 3-body | 7,463,678 | 125,380,883,295 | $120 \times 6$ |
| $^{10}$He | 8 | 3-body | 15,631,854 | 247,561,956,281 | $276 \times 6$ |
| $^8$Li | 8 | 3-body | 17,684,028 | 399,084,871,767 | $496 \times 6$ |

## 5 Numerical Results and Computational Performance

In the following, we present computational examples performed on Edison, a Cray XC30 system with 5576 computational nodes, at the National Energy Research Scientific Computing Center (NERSC).[6] Each computational node has 24 cores, located on two twelve-core Intel "Ivy Bridge" 2.4 GHz processors, and has 64 GB DDR3 1866 MHz memory. Each physical core has its own L1 and L2 caches, with 64 KB (32 KB instruction cache, 32 KB data) and 256 KB, respectively, and can run one or two user threads (i.e., hyperthreading). A 30 MB L3 cache is shared by twelve cores on the "Ivy Bridge" processor. The computational nodes on Edison are connected by Cray Aries network with Dragonfly topology with 23.7 TB/s global bandwidth. The MFDn code is compiled with Intel Compiler version 15.0.1, and linked with the Cray Scientific Libraries package (LibSci) version 13.3.0 and Cray's MPI library version 7.3.1. We run the code using four MPI processes per node, and six OpenMP threads per MPI process—no hyperthreading is utilized.

We run MFDn on eight test cases as shown in Table 1, using both Lanczos and LOBPCG solvers. As indicated in the third column of Table 1, we select four test cases with the two-body nucleon–nucleon interaction alone and four test cases that include a three-body interaction. Each of the four cases is paired in computational intensity as indicated by the number of nonzero elements in the lower triangular portion of $\hat{H}$. The eight lowest eigenvalues and the associated eigenvectors are calculated for each test case. The convergence criterion for the eigensolvers is that the relative residual norm drops below $10^{-6}$.

### 5.1 Impact of preconditioning and initial guess in LOBPCG

To demonstrate how the techniques discussed in Section 3 help improve the efficiency of LOBPCG, we perform several runs for the $^7$Li example with different settings: with/without preconditioning, and with a good initial guess from a smaller configuration space or with a random initial guess. When adopting a good initial guess from a smaller configuration space, we run LOBPCG with $N_{\max} = 8$ and $N_{\max} = 10$, and use the results as initial guesses for the LOBPCG method $N_{\max} = 10$ and $N_{\max} = 12$, respectively. The iteration counts (for $N_{\max} = 12$) and execution times of the diagonalization phase are summarized in Table 2.

Both preconditioning and a good initial guess can effectively reduce the number of iterations and the execution times for LOBPCG. If neither technique is adopted, LOBPCG becomes slower than Lanczos for this example. Hence, an efficient SpMM implementation alone is not sufficient for an efficient block eigensolver. Note that the execution time for the LOBPCG runs with a good initial guess includes the time for constructing the initial guess. But the iteration count refers to the LOBPCG run for the $N_{\max} = 12$ case. In the run with preconditioning and good initial guess, it takes 38.2 seconds, about 18% of the overall diagonalization time, to perform two LOBPCG runs on smaller configuration spaces. The benefit of the runs in smaller configuration spaces is the reduced iteration counts and thereby reduced execution times when using the LOBPCG solver for the largest configuration space with $N_{\max} = 12$.

---

[6]http://www.nersc.gov/systems/edison-cray-xc30/

Table 2: Execution times (in seconds) for different settings of LOBPCG for the $^7$Li example at $N_{\max} = 12$ using 496 MPI processes. Execution times cover the overall diagonalization time and the IO (i.e., writing the solution to files) time. The timing from the Lanczos algorithm is also listed for comparison.

| LOBPCG | | iterations | time |
|---|---|---|---|
| preconditioning | good initial guess | | |
| No | No | 154 | 447.21 |
| No | Yes | 81 | 250.75 |
| Yes | No | 105 | 304.74 |
| Yes | Yes | 53 | 206.67 |
| Lanczos | | 400 | 387.09 |

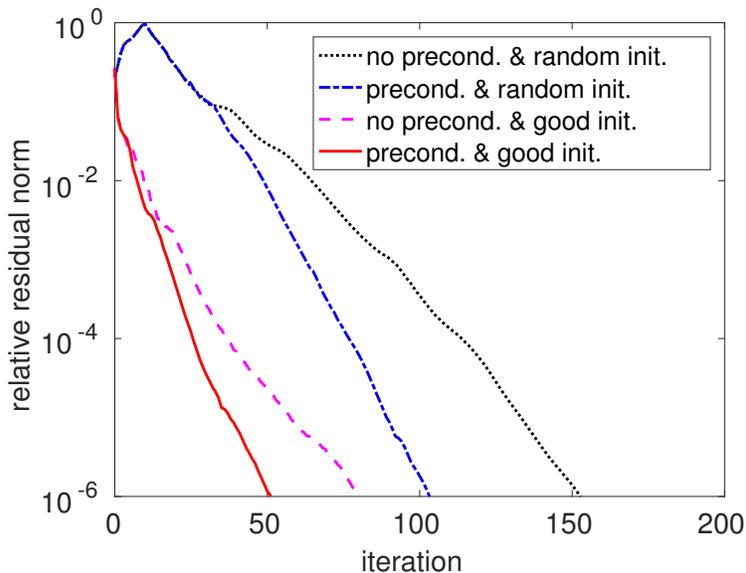

Figure 6: Convergence history of LOBPCG for the $^7$Li example. Both a good initial guess and preconditioning largely improve the convergence rate of LOBPCG.

Fig. 6 shows the convergence history of the LOBPCG algorithm for these different settings. It can be clearly seen from the figure that the use of a random initial guess results in taking more than 25 iterations just to bring the residual norm down to the same level satisfied by a good initial guess. Also, the slopes of the curves with preconditioning are steeper than those without preconditioning, indicating that preconditioning improves the convergence rate. Hence both techniques are helpful for reducing the total number of iterations. Finally, we remark that other examples have similar behaviors using LOBPCG vs. Lanczos solvers, therefore we present these detailed convergence patterns only for this $^7$Li example.

Table 3: Execution times of diagonalization phase (in seconds). Three different strategies of choosing initial guess in LOBPCG are tested: Column 3 ("smaller $N_{\max}$")—use solution from a smaller configuration space as the initial guess; Column 4 ("$\hbar\omega' = 0.8\hbar\omega$")—use solution from a different basis parameter $\hbar\omega$ as the initial guess; Column 5 ("different $V$")—use solution with two-body potential as the initial guess in a run with three-body potential. The execution times of Lanczos algorithm are also listed for comparison. Speedup are measured as time ratio between Lanczos and LOBPCG.

| Nucleus | Lanczos $n_{\text{it}}$[7] | Lanczos time | smaller $N_{\max}$ $n_{\text{it}}$ | smaller $N_{\max}$ time | $\hbar\omega' = 0.8\hbar\omega$ $n_{\text{it}}$ | $\hbar\omega' = 0.8\hbar\omega$ time | different $V$ $n_{\text{it}}$ | different $V$ time | speedup |
|---|---|---|---|---|---|---|---|---|---|
| $^9$Li | 240 | 209.8 | 36 | 129.4 | 37 | 124.3 | 47 | 150.1 | $1.4 \sim 1.7$ |
| $^{10}$B | 280 | 344.9 | 35 | 187.3 | 48 | 223.0 | 40 | 203.3 | $1.6 \sim 1.9$ |
| $^8$Be | 480 | 531.3 | 53 | 194.1 | 59 | 192.2 | 62 | 203.1 | $2.6 \sim 2.8$ |
| $^7$Li | 400 | 389.7 | 56 | 201.3 | 60 | 213.7 | 60 | 177.7 | $1.8 \sim 2.2$ |
| $^6$Li | 240 | 156.1 | 33 | 46.5 | 34 | 58.7 | 30 | 51.3 | $2.7 \sim 3.3$ |
| $^8$He | 320 | 247.7 | 66 | 111.4 | 59 | 93.8 | 56 | 88.2 | $2.2 \sim 2.8$ |
| $^{10}$He | 480 | 319.4 | 122 | 169.2 | 115 | 156.4 | 126 | 174.2 | $1.8 \sim 2.0$ |
| $^8$Li | 280 | 169.3 | 38 | 58.7 | 37 | 53.7 | 35 | 56.2 | $2.9 \sim 3.2$ |

## 5.2 Comparison between LOBPCG and Lanczos

In the following we compare the LOBPCG method (using preconditioning and a good initial guess from smaller configuration spaces) with the Lanczos method. The execution times as well as the iteration counts for all test cases are shown in Table 3. Though LOBPCG performs much more computation in SpMM (roughly speaking, each LOBPCG iteration requires applying SpMM to twelve vectors as we use 50% more trial vectors compared to the number of desired eigenvectors to enhance the robustness), it consistently outperforms the Lanczos algorithm: up to and beyond $3\times$ improvements in the wallclock execution time is observed. (In terms of flop rate, the performance improvement is even larger: with LOBPCG we perform more flops and consume less time than Lanczos.)

In Table 3 we also show execution times for other choices for the initial guess. Since we often perform a scan of several basis parameters $\hbar\omega$, it is straightforward to use solutions from a previous run as an initial guess for the next run. In Table 3, we show numbers of iterations and times for the case where the preconditioner is adopted from results at a basis parameter (called $\hbar\omega'$) that is 80% of the current value of basis parameter. This generally shows a similar reduction in execution time as using solutions from a smaller configuration space. Using a solution obtained with a different nuclear potential also shows similar improvements: in Table 3 we show in the column labeled "different $V$" the iteration count and execution time using solution from a different NN potential (top four entries), and from the NN-only part of an NN-plus-3N potential (bottom four entries). This option can be especially useful when fitting parameters in the NN and/or 3N potential. (Of course, the closer the initial guess is to the actual desired solutions the better.)

To analyze why LOBPCG outperforms Lanczos in all test cases, we illustrate in Fig. 7 a detailed profiling on execution times. Overall, SpMV (SpMM) dominates the execution time for Lanczos (LOBPCG), especially in runs involving three-body interactions (results in the bottom row of Fig. 7). The use of SpMM in LOBPCG is clearly more efficient than the SpMV operations performed in Lanczos as SpMM performs more computational work but consumes less execution time. For runs with two-body interactions (top row of panels in Fig. 7), the orthogonalization cost is also non-negligible. This is likely due to the relatively large amount of flops required in the orthogonalization process compared to three-body runs. The communication cost in orthogonalization for the Lanczos method is sometimes quite significant. As a comparison, orthogonalization in LOBPCG is performed using level-3 operations, hence uses more efficient communication. As desired, the preconditioning phase in LOBPCG is cost effective, especially for large systems.

---

[7]The convergence of the Lanczos algorithm is checked every 10 iterations at the beginning, and then every 20 iterations after 100 iterations, and then every 40 iterations after 200 iterations, etc.

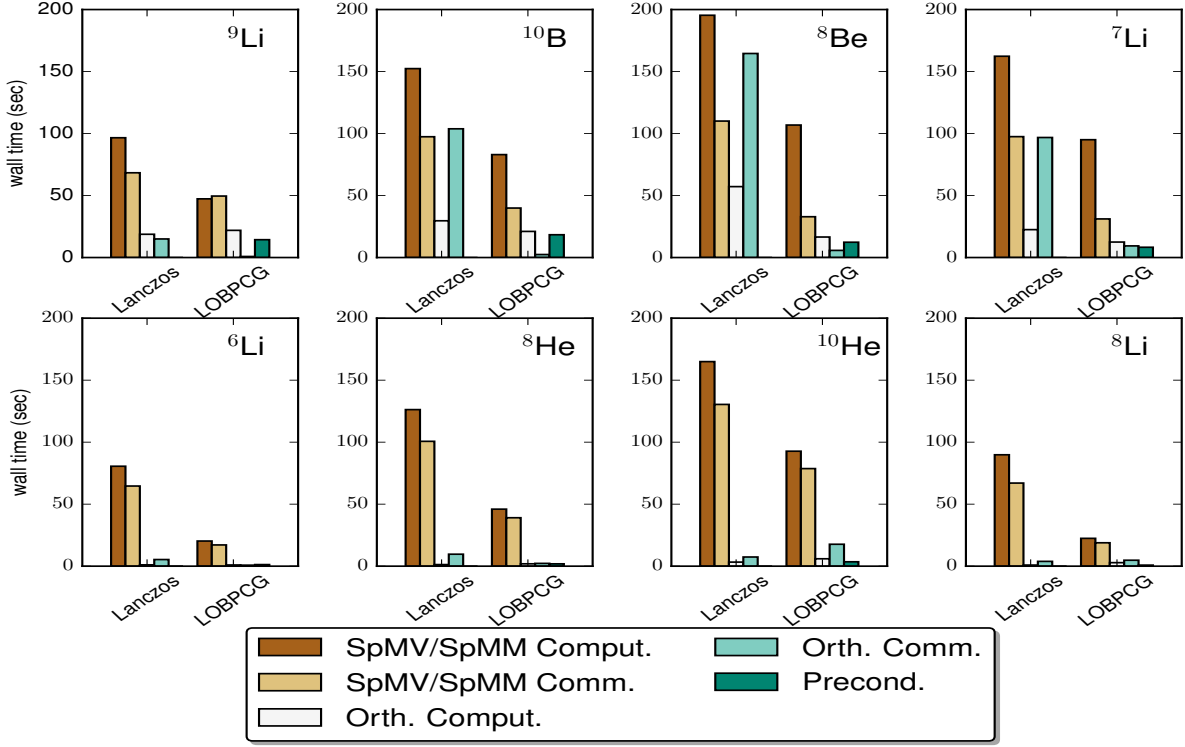

Figure 7: Detailed profiling on execution times for Lanczos and LOBPCG (using initial guess from a smaller $N_{\max}$). LOBPCG outperforms Lanczos in both communication and computations. In LOBPCG, the time consumed by preconditioning is relatively small compared to other components.

## 6 Conclusions

We described recently developed techniques for improving the performance of nuclear configuration interaction calculation. These techniques include

1. using a block iterative method such as the LOBPCG algorithm that exposes a higher level of concurrency, and can take advantage of the memory hierarchy on modern microprocessors to increase arithmetic intensity relative to data movement;

2. choosing an appropriate starting guess of the eigenvectors;

3. constructing an effective preconditioner.

The latter two techniques exploit special structures of the nuclear many-body Hamiltonian which we described in detail. Understanding these structures is also key to developing an effective data distribution and parallelization strategy to achieve high performance on massively parallel multi-core supercomputers. The performance improvements achieved by using our developments are demonstrated by several examples.

## Acknowledgments

This work was supported in part by the U.S. Department of Energy (DOE) under grants No. DESC0008485 (SciDAC/NUCLEI) and DE-FG02-87ER40371. Computational resources were provided by the National Energy Research Scientific Computing Center (NERSC), which is supported by the U. S. Department of Energy under Contract No. DE-AC02-05CH11231. H. M. Aktulga was also supported by a research grant from Michigan State University.